\documentclass[reqno, 11pt]{amsart}

\usepackage[utf8]{inputenc}
\usepackage{amssymb}
\usepackage{mathrsfs}
\usepackage{graphicx}
\usepackage{latexsym}
\usepackage{stmaryrd}
\usepackage{enumitem}
\usepackage[bookmarks,hyperfootnotes=false]{hyperref}
\usepackage{multirow}
\usepackage{xcolor}
\usepackage{caption}
\colorlet{darkishRed}{red!80!black}
\colorlet{darkishBlue}{blue!60!black}
\colorlet{darkishGreen}{green!60!black}
\hypersetup{
    draft = false,
    bookmarksopen=true,
    colorlinks,
    linkcolor={red!60!black},
    citecolor={green!60!black},
    urlcolor={blue!60!black}
}
\usepackage[nameinlink, capitalise, noabbrev]{cleveref}
\crefformat{enumi}{#2#1#3}
\crefformat{equation}{#2(#1)#3}
\usepackage[msc-links]{amsrefs} 
\usepackage{doi}

\renewcommand{\PrintDOI}[1]{\doi{#1}}

\let\setminus=\smallsetminus
\usepackage{tikz}
\usepackage{tikz-cd}
\usetikzlibrary{calc,through,intersections,arrows, trees, positioning, decorations.pathmorphing, cd, shapes.geometric}
\usepackage{relsize}
\usepackage{comment}
\usepackage{svg}
\usepackage{mathtools}
\usepackage{nccmath}
\usepackage{pifont}
\usepackage{standalone}

\usepackage{geometry}
\let\setminus=\smallsetminus
\renewcommand{\leq}{\leqslant}

\renewcommand{\le}{\leq}

\newcommand{\modG}{\vert G\vert}
\newcommand{\modGn}{\vert G_n\vert}

\renewcommand{\subset}{\subseteq}
\renewcommand{\supset}{\supseteq}

\newcommand{ \N } { \mathbb{N} }

\makeatletter

\def\calCommandfactory#1{%
   \expandafter\def\csname c#1\endcsname{\mathcal{#1}}}
\def\frakCommandfactory#1{%
   \expandafter\def\csname frak#1\endcsname{\mathfrak{#1}}}

\newcounter{ctr}
\loop
  \stepcounter{ctr}
  \edef\X{\@Alph\c@ctr}
  \expandafter\calCommandfactory\X
  \expandafter\frakCommandfactory\X
  \edef\Y{\@alph\c@ctr}
  \expandafter\frakCommandfactory\Y
\ifnum\thectr<26
\repeat

\setenumerate{label={\normalfont (\roman*)}}

\usepackage[presets={vec-cev,abc,ABC,cAcBcC}]{letterswitharrows}

\newtheorem{theorem}{Theorem}[section] 
\newtheorem{corollary}[theorem]{Corollary}
\newtheorem{lemma}[theorem]{Lemma}

\newtheorem*{question*}{Question}
\newtheorem{mainresult}{Theorem}

\theoremstyle{definition}
\newtheorem{example}[theorem]{Example}

\theoremstyle{remark}

\newtheorem*{ack}{Acknowledgement}

\usepackage[nameinlink, capitalise, noabbrev]{cleveref}
\crefformat{enumi}{#2(#1)#3}

\title[The Lov\'{a}sz-Cherkassky theorem in $\left|G\right|$]{The Lov\'{a}sz-Cherkassky theorem for\\locally finite graphs with ends\\\medskip -- Note --}

\author{Raphael W. Jacobs}
\thanks{The first author gratefully acknowledges support by doctoral scholarships of the Studienstiftung des deutschen Volkes and the Cusanuswerk -- Bisch\"{o}fliche Studienf\"{o}rderung.}

\author{Attila Jo\'{o}}
\thanks{The second author would like to thank the generous support of the Deutsche Forschungsgemeinschaft (DFG, German
Research Foundation)-513023562, Alexander 
von Humboldt Foundation and NKFIH OTKA-129211}
\address{Attila Jo\'{o}:
Universität Hamburg, Department of Mathematics, Bundesstra{\ss}e 55 (Geomatikum), 20146 Hamburg, Germany}
\email{attila.joo@uni-hamburg.de}

\author{Paul Knappe}
\address{Raphael W. Jacobs, Paul Knappe, Jan Kurkofka, Ruben Melcher: Universität Hamburg, Department of Mathematics, Bundesstra{\ss}e 55 (Geomatikum), 20146 Hamburg, Germany}
\email{raphael.jacobs@uni-hamburg.de, paul.knappe@uni-hamburg.de}
\email{j.lastname@bham.ac.uk, ruben.melcher@uni-hamburg.de}

\author{Jan Kurkofka}

\author{Ruben Melcher}

\keywords{Lov\'{a}sz-Cherkassky theorem, infinite graph, Freudenthal compactification, edge-connectivity}
\subjclass[2020]{Primary: 05C63, 05C40, 05C38.} 

\begin{document}

\begin{abstract}
    Lov\'{a}sz and Cherkassky discovered independently that, if $G$ is a finite graph and $T\subseteq V(G)$ such that the degree $d_G(v)$ is even for every vertex $v\in V(G)\setminus T$, then the maximum number of edge-disjoint paths which are internally disjoint from~$T$ and connect distinct vertices of $T$ is equal to $\frac{1}{2} \sum_{t\in T}\lambda_G(t, T\setminus \{t\})$ (where $\lambda_G(t, T\setminus \{t\})$ is the size of a smallest cut that separates $t$ and $T\setminus\{t\}$). From another perspective, this means that for every vertex $t\in T$, in any optimal path-system there are $\lambda_G(t, T\setminus \{t\})$ many paths between $t$ and~$T\setminus\{t\}$. 
    We extend the theorem of Lov\'{a}sz and Cherkassky based on this reformulation to all locally-finite infinite graphs and their ends. 
    In our generalisation, $T$ may contain not just vertices but ends as well, and paths are one-way (two-way) infinite when they establish a vertex-end (end-end) connection.
\end{abstract}

\vspace*{-2cm}
\maketitle

\vspace*{-.5cm}
\section{Introduction}

\noindent 
A non-trivial path $P$ is a $T$-\emph{path} for a set $T$ of vertices if $P$ has its endvertices but no inner vertex in~$T$.
For disjoint vertex sets $X$ and~$Y$ in a graph~$G$, we write $\lambda_G(X,Y)$ for the size of a smallest cut in~$G$ that separates~$X$ and~$Y$.

Now let $T$ be any set of vertices in a finite graph~$G$.
In a set $\mathcal{P}$ of edge-disjoint $T$-paths, there are at most $\lambda_G(t, T\setminus \{t\})$ many paths that link a vertex~$t\in T$ to~$T\setminus\{t\}$. It follows that 
$\left|\mathcal{P}\right|\leq \frac{1}{2} \sum_{t\in T}\lambda_G(t, T\setminus \{t\})$. The question about the sharpness of this upper bound can be formulated in the following structural way.
For every vertex $t\in T$, let $\cP_t$ be a set of $\lambda_G(t, T\setminus \{t\})$ many edge-disjoint $t$--$(T\setminus\{t\})$ paths.

\begin{question*}
    Can we choose the paths in the sets $\cP_t$ for each vertex $t\in T$ in such a way that the union of all the sets $\cP_t$ is an edge-disjoint path-system? 
\end{question*}

\noindent
Clearly, the answer is no: the three leaves of a star $K_{1,3}$ form a set~$T$ where each set $\cP_t$ must consist of a single path, but the union $\bigcup_{t\in T}\cP_t$ always contains two distinct paths that share an edge, no matter how we choose the paths in each set~$\cP_t$.
Lov\'{a}sz and Cherkassky independently showed that, perhaps surprisingly, the answer is yes under the additional assumption that the graph~$G$ is \emph{inner-Eulerian} for~$T$ in that every vertex of~$G$ which is not in~$T$ has even degree in~$G$.

\begin{theorem}[Lov\'{a}sz-Cherkassky Theorem \cites{lovasz1976some, cherkassky1977asol}] \label{t:LCh}
    Let $G$ be any finite graph, and let $T\subseteq V(G)$ such that $G$ is inner-Eulerian for~$T$.
    Then the maximum number of pairwise edge-disjoint $T$-paths in~$G$ is equal to
	\[\frac{1}{2}\sum_{t\in T}\lambda_G(t, T\setminus\{t\}).\]
\end{theorem}

Diestel asked whether \cref{t:LCh} extends to all locally-finite infinite graphs and their ends~\cite{DiestelPersComm}.
In this note, we answer his question in the affirmative.
For this, we employ the Freudenthal compactification in the spirit of \cites{diestel2011topapr1, diestel2010topapr2, diestel2012topapr3}, as customary in the study of locally-finite infinite graphs~\cites{diestel2004topological,diestel2004infinitecyc1,diestel2004infinitecyc2,bruhn2005menger,bruhn2007end,diestel2003countable}.
For infinite graphs and their ends, we follow and assume familiarity with the terminology in~\cite{diestelbook}*{§8}, in particular in~§8.6.

We allow a graph to have parallel edges, but we do not allow any loops; if a graph has no parallel edges, we call it \emph{simple}.
The \emph{degree} of a vertex~$v$ in a graph~$G$ is the number $d_G(v)\in\N\cup\{\infty\}$ of edges of~$G$ incident with~$v$.
If all the vertices of $G$ have finite degree, then we say that $G$ is \emph{locally finite}.
Note that in a locally finite graph, there can be only finitely many parallel edges between any two vertices.
We write $\hat{V}(G)$ for the union of the vertex set $V(G)$ of~$G$ and the set $\Omega(G)$ of all ends of~$G$.
An arc~$A$ in the end compactification~$\modG$ of a locally-finite connected graph~$G$ is a $T$-\emph{arc} for a set $T\subseteq \hat{V}(G)$ if $A$ has its endpoints but no inner points in~$T$.
An arc~$A$ is an $X$--$Y$ \emph{arc} between two sets $X$~and~$Y$ if $A$ intersects $X$ precisely in one endpoint and $Y$ precisely in the other. 
We call an arc \emph{graphic} if it is defined by a finite graph-theoretic path, a ray, or a double ray with its tails in distinct ends.
Two arcs in~$\modG$ are \emph{edge-disjoint} if they do not meet in inner points of edges.

A subset $X\subset \hat{V}(G)$ \emph{lives in} a subgraph $C\subset G$ or a vertex set $C\subset V(G)$ if all the vertices of~$X$ lie in~$C$ and all the rays of ends in~$X$ have tails in~$C$ or~$G[C]$, respectively.
A finite cut $F$ of a graph $G$ is an $X$--$Y$ \emph{cut} for two sets $X,Y\subset\hat{V}(G)$ if $X$ and $Y$ live in distinct sides of~$F$.
If $Y\subset \hat{V}(G)$ is a set of vertices and ends of a locally finite graph~$G$, and $x\in\hat{V}(G)$ is not contained in the closure of~$Y$ in~$\modG$, then $G$ admits a (finite) $x$--$Y$ cut, and we denote the least size of such a cut by~$\lambda_G(x,Y)$.

Our main result reads as follows:
\begin{mainresult} \label{t:MainResult}
    Let $G$ be any locally-finite graph, and let $T\subseteq \hat{V}(G)$ be discrete in $\modG$. If every finite cut of $G$ such that $T$ lives on one of its sides is even, then $\modG$ contains a set $\cA$ of pairwise edge-disjoint graphic $T$-arcs such that for every $t \in T$, the number of $t$--$(T\setminus\{t\})$ arcs in~$\cA$ is equal to~$\lambda_G(t,T\setminus\{t\})$.
\end{mainresult}

\noindent On first sight, one might wonder why the notion of `inner-Eulerian' from the premise of \cref{t:LCh} has been replaced in the premise of \cref{t:MainResult} with a new condition on the vertex sets~$X$.
In short, in finite graphs `inner-Eulerian' implies the new premise, so \cref{t:MainResult} is more general; see \cref{app:degrees} for a discussion of the necessity of the new condition for end-compactifications of locally finite graphs.

The assumption that~$T$ is a discrete subset of~$\modG$ naturally arises here as it is equivalent to asking that there exists a $t$--$(T\setminus\{t\})$ cut in~$G$ for each end~$t\in T$, which precisely ensures that $\lambda_G(t,T\setminus\{t\})$ is defined for all~$t\in T$.
We remark that the assumption that~$T$ is discrete is also motivated by the work of Bruhn, Diestel, and Stein~\cite{bruhn2005menger}, which is a generalisation of the Erd\H{o}s-Menger theorem by Aharoni and Berger~\cite{aharoni2009menger} from infinite graphs to infinite graphs and their ends, under a similar assumption on the ends which implies discreteness in our setting.
(Diestel discusses this assumption in detail in~\cite{diestel2003countable}*{§3}.) 

We conclude the introduction with an example that discusses why a natural weakening of the discreteness-assumption in \cref{t:MainResult} cannot be made.

\begin{figure}[ht]
	\centering
	\begin{tikzpicture}[>=stealth',every node/.style={scale=0.62}]
	
	    \node (u4r) at (5,0) {};
	    \node (v4r) at (5,-1) {};
	    \node (u0l) at (-1,0) {};
	    \node (v0l) at (-1,-1) {};
		
		\node[circle, fill=black] (u0) at (0,0) {};
		\node[circle, fill=black] (u1) at (1,0) {};
		\node[circle, fill=black] (u2) at (2,0) {};
		\node[circle, fill=black] (u3) at (3,0) {};
		\node[circle, fill=black] (u4) at (4,0) {};
		
		\node[circle, fill=black] (endr) at (6,-.5) {};
		\node[circle, fill=black] (endl) at (-2,-.5) {};
		
		\node[circle, draw] (v0) at (0,-1) {};
		\node[circle, draw] (v1) at (1,-1) {};
		\node[circle, draw] (v2) at (2,-1) {};
		\node[circle, draw] (v3) at (3,-1) {};
		\node[circle, draw] (v4) at (4,-1) {};
		
		\draw[->] (u4) edge (u4r);
		\draw[->] (v4) edge (v4r);
		\draw[->] (u0) edge (u0l);
		\draw[->] (v0) edge (v0l);
		
		\draw  (v0) edge (v1);
		\draw  (v1) edge (v2);
		\draw  (v2) edge (v3);
		\draw  (v3) edge (v4);
		\draw  (u0) edge (u1);
		\draw  (u1) edge (u2);
		\draw  (u2) edge (u3);
		\draw  (u3) edge (u4);
		
		\draw (u0) edge [bend left] (v0);
		\draw (u0) edge [bend right] (v0);
		\draw (u1) edge [bend left] (v1);
		\draw (u1) edge [bend right] (v1);
		\draw (u2) edge [bend left] (v2);
		\draw (u2) edge [bend right] (v2);
		\draw (u3) edge [bend left] (v3);
		\draw (u3) edge [bend right] (v3);
		\draw (u4) edge [bend left] (v4);
		\draw (u4) edge [bend right] (v4);
	\end{tikzpicture}
	\caption{$T$ consists of the black vertices and ends} \label{fig:FailureOne}
\end{figure}
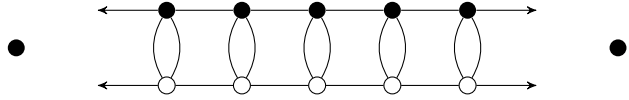 

\begin{example}
    We claim that it is not possible to drop in \cref{t:MainResult} the requirement that $T$ is discrete in~$\modG$ and, instead, replace $\lambda_G(t,T\setminus\{t\})$ in the wording of the theorem with the maximum number $\mu_G(t,T\setminus\{t\})$ of pairwise edge-disjoint graphic $t$--$(T\setminus\{t\})$ arcs in~$\modG$.
    Indeed, let $G$ be obtained from the double ladder by duplicating each rung, and let $T$ consist of both ends of~$G$ together with the vertices of one of the two main double rays; see~\cref{fig:FailureOne}.
    Then
    \begin{equation*}
        \mu_G(t, T\setminus \{t\}) = 
        \begin{cases}
            4 & \mbox{if } t\in T\cap V(G), \\
            1 & \mbox{if } t\in T\cap \Omega(G).
        \end{cases}
    \end{equation*}
    But any $t$--$(T\setminus\{t\})$ arc for some $t \in T \cap \Omega(G)$ does already preclude the existence of four $t'$--$(T\setminus\{t'\})$ arcs for some of the $t' \in T \cap V(G)$, and hence there is no desired arc-system for $T$.
\end{example}

\subsection{The Lov\'{a}sz-Cherkassky Theorem for infinite graphs without ends}

The second author has extended \cref{t:LCh} to infinite -- not necessarily locally finite -- graphs where the set~$T$ is a countably infinite set of vertices, but does not contain ends.
A cut $F$ is said to lie \emph{on} a set $\cP$ of edge-disjoint paths in~$G$ if $F$ consists of a choice of exactly one edge from each path in~$\cP$.
For a vertex set $X\subset V(G)$, we denote by~$d_G(X)$ the number of edges of~$G$ between $X$ and its complement $V(G)\setminus X$.
Note that this notation is consistent with the above definition of the degree $d_G(v)$ of a vertex $v$ of $G$ in that $d_G (v) = d_G(\{v\})$ for every vertex $v \in G$.

\begin{theorem}[\cite{lovcherinf}*{Theorem 1.3}] \label{t:LChInfinite}
	Let $G$ be any graph, and let $T\subseteq V(G)$ be a countable vertex set such that there is no $X\subseteq V(G)\setminus T$ 
	for which $d_G(X)$ is an odd natural number. 
	Then $G$ contains a set	$\cP$ of edge-disjoint $T$-paths such that for each vertex $t\in T$, the graph $G$ contains a $t$--$(T\setminus\{t\})$ cut on the set of $t$--$(T\setminus\{t\})$ paths in~$\cP$.
\end{theorem}

\noindent We will use \cref{t:LChInfinite} in the proof of \cref{t:MainResult}.

\cref{t:LChInfinite} compares to \cref{t:MainResult} as follows.
On the one hand, \cref{t:LChInfinite} can be applied in the setting of countable graphs that are not locally finite, where \cref{t:MainResult} cannot be applied.
In the setting of locally finite graphs, on the other hand, \cref{t:MainResult} is more general than \cref{t:LChInfinite}, as it allows $T$ to consist of vertices and ends alike.
Indeed, the locally-finite version of \cref{t:LChInfinite} can easily be re-obtained from \cref{t:MainResult} since every set of vertices of a locally finite graph~$G$ is discrete in~$|G|$ and does not contain any other vertices in its closure.

\section{Proof of the main result} \label{sec:Proof}

\noindent We need one auxiliary result for the proof \cref{t:MainResult}, and to state this lemma we make the following definition. 
If $S$ is a set of vertices of a graph~$G$ and $\omega$ is an end of~$G$, then by an $S$--$\omega$ \emph{ray} we mean a ray which has precisely its first vertex in~$S$ and belongs to~$\omega$.

\begin{lemma}[{\cite{bruhn2007end}*{Lemma 10}}] \label{lemma:MengerVerticesEnd}
    Let $G$ be a locally finite connected graph, let $\omega$ be an end of $G$, and let $S$ be a finite set of vertices in $G$.
    Then the maximum number of edge-disjoint $S$--$\omega$ rays is equal to the minimum size of a cut that separates $S$ and~$\omega$.
\end{lemma}

\noindent Originally, \cref{lemma:MengerVerticesEnd} has been proved only for simple graphs.
However, it extends to graphs with parallel edges: Given a graph with parallel edges, just subdivide each edge once, apply the original result, and then suppress all subdividing vertices.

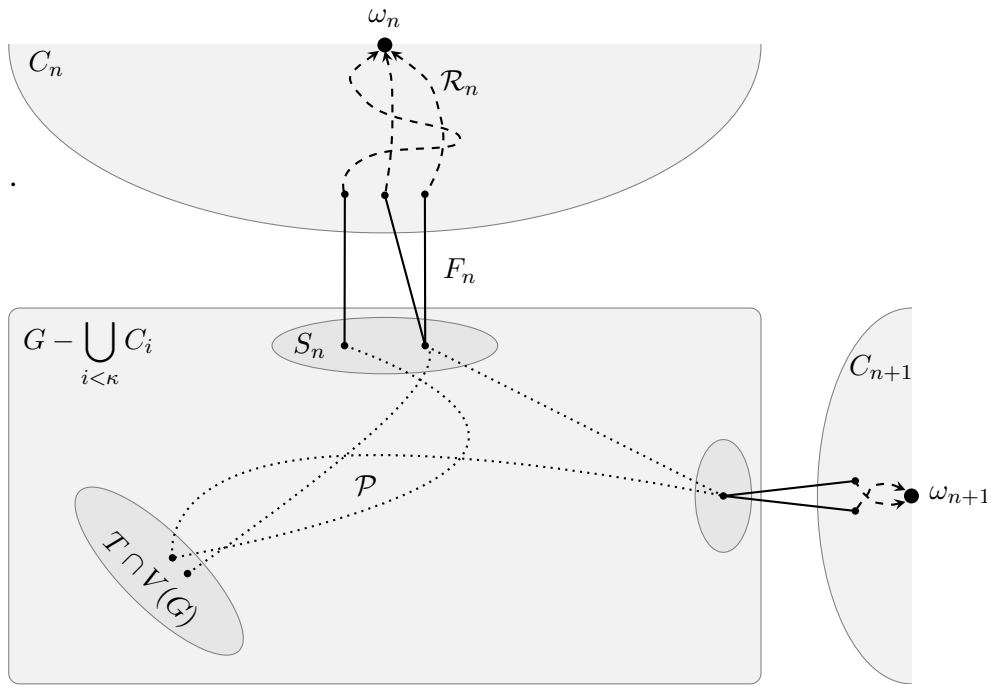
\begin{figure}[t]
    \centering
    \begin{tikzpicture}
	\clip (-6,6.5) rectangle (8, -4);

    \edef \x{6};

    \node[outer sep = -0.5mm] (aux1) at (0,\x) {};
	\node[shape = coordinate] (center) at (0,0) {};
	
	\begin{scope}[yshift=\x cm]
		\clip(-5,0)rectangle(5,-3);
		\node[ellipse,draw=black!50, minimum width = 10cm, minimum height = 5cm, outer sep = -0.5cm, fill = black!05] (Cninner) at (aux1){};
	\end{scope}
	
	\node [label = above: $\omega_n$] at (aux1)[fill, circle, inner sep = 2pt]{};
	
	\node[shape = coordinate] (s0) at (Cninner.255) {};
	\node[shape = coordinate] (s1) at (Cninner.270) {};
	\node[shape = coordinate] (s2) at (Cninner.285) {};
	
	\node[rectangle, rounded corners, draw=black!50, minimum width = 10cm, minimum height = 5cm, outer sep = -0.5cm, fill = black!05] (GLeftover) at (center){};
	
	\node[shape = coordinate] (r0) at (GLeftover.105) {};
	\node[shape = coordinate] (r1) at (GLeftover.90) {};
	\node[shape = coordinate] (r2) at (GLeftover.75) {};

	\node[ellipse,draw=black!50, minimum width = 3cm, minimum height = 0.75cm, fill = black!10] (Sn) at (r1){};

	\draw[thick] plot [smooth,tension=1] coordinates {(s0) (r0)};
	\draw[thick] plot [smooth,tension=1] coordinates {(s1) (r2)};
	\draw[thick] plot [smooth,tension=1] coordinates {(s2) (r2)};
	
	\draw[<-,>=stealth,thick,dashed] plot [smooth, tension=1] coordinates {(aux1.225) (-0.25,5.25) (1,4.75) (-0.25, 4.5) (s0)};
	\draw[<-,>=stealth,thick,dashed] plot [smooth, tension=1] coordinates {(aux1.270) (0.1,5) (s1)};
	\draw[<-,>=stealth,thick,dashed] plot [smooth, tension=1] coordinates {(aux1.315) (0.75,5) (s2)};
	
	\node at (1,3) {$F_n$};   
	\node at (-4.5,5.75) {$C_n$};
	\node at (1,5.5) {$\mathcal{R}_n$};
	\node at (-1,2) {$S_n$};

	\edef \y{7};
	\node[outer sep = -0.5mm] (aux2) at (\y,0) {};
	
	\begin{scope}[xshift=\y cm]
		\clip(0,5)rectangle(-2.5,-2.5);
		\node[ellipse,draw=black!50, minimum width = 2.5cm, minimum height = 5cm, outer sep = -0.5cm, fill = black!05] (Cn1inner) at (aux2){};
	\end{scope}
	
	\node at (6.6,1.7) {$C_{n+1}$};
	
	\node [label = right: $\omega_{n+1}$] at (aux2)[fill, circle, inner sep = 2pt]{};
	
	\node[shape = coordinate] (s10) at (Cn1inner.195) {};
	\node[shape = coordinate] (s11) at (Cn1inner.-195) {};

	\node[shape = coordinate] (r10) at (GLeftover.0) {};

	\node[ellipse,draw=black!50, minimum width = 0.75cm, minimum height = 1.5cm, fill = black!10] (Sn) at (r10){};

	\draw[thick] plot [smooth,tension=1] coordinates {(s10) (r10)};
	\draw[thick] plot [smooth,tension=1] coordinates {(s11) (r10)};

	\draw[<-,>=stealth,thick,dashed] plot [smooth, tension=1] coordinates {(aux2.135) (6.55,0.15) (6.4,0) (s10)};
	\draw[<-,>=stealth,thick,dashed] plot [smooth, tension=1] coordinates {(aux2.225) (6.6,-0.1) (6.4,0) (s11)};
	
	\node[ellipse, minimum width = 15cm, minimum height = 11cm] (circ) at (center){};
	\node(dot1) at (circ.140)[fill, circle, inner sep = 0.5pt]{};
	\node(dot2) at (circ.145)[fill, circle, inner sep = 0.5pt]{};
	\node(dot3) at (circ.150)[fill, circle, inner sep = 0.5pt]{};
	
	\node(dot4) at (circ.335)[fill, circle, inner sep = 0.5pt]{};
	\node(dot5) at (circ.330)[fill, circle, inner sep = 0.5pt]{};
	\node(dot6) at (circ.325)[fill, circle, inner sep = 0.5pt]{};
	
	\node[ellipse, draw = black!50, rotate= 135, minimum width = 3cm, minimum height = 1cm, fill = black!10, outer sep = -0.25cm] (T) at (-3,-1){};
	\node[rotate = 315] at (-3.125,-1.125) {$T \cap V(G)$};
	
	\node(t0) at (T.270){};
	\node(t1) at (T.220){};
	
	\draw[thick,dotted] plot [smooth, tension=1] coordinates {(r0) (1,0.5) (t0)};
	\draw[thick,dotted] plot [smooth, tension=1] coordinates {(r2) (0,1)(t1)};
	\draw[thick,dotted] plot [smooth, tension=1] coordinates {(r2) (r10)};
	\draw[thick,dotted] plot [smooth, tension=1] coordinates {(t0) (-1,0.5) (r10)};
	
	\node [] at (s0)[fill, circle, inner sep = 1pt]{};
	\node [] at (s1)[fill, circle, inner sep = 1pt]{};
	\node [] at (s2)[fill, circle, inner sep = 1pt]{};
	\node [] at (r0)[fill, circle, inner sep = 1pt]{};
	\node [] at (r2)[fill, circle, inner sep = 1pt]{};
	\node [] at (s10)[fill, circle, inner sep = 1pt]{};
	\node [] at (s11)[fill, circle, inner sep = 1pt]{};
	\node [] at (r10)[fill, circle, inner sep = 1pt]{};
	\node [] at (t0)[fill, circle, inner sep = 1pt]{};
	\node [] at (t1)[fill, circle, inner sep = 1pt]{};
	
	\node at (-0.25,0.15) {$\mathcal{P}$};
	\node at (-3.94,1.9) {$G - \displaystyle\bigcup_{i < \kappa} C_i$};
\end{tikzpicture}
    \caption{Situation in the proof of the main result}
    \label{fig:MainProof}
\end{figure}

\begin{proof}[Proof of \cref{t:MainResult}]
    Without loss of generality, we may assume that $G$ is connected.
    
	Let us first show that $T$ is countable.
	Since~$G$ is locally finite, $\modG$~is second-countable, meaning that the topology on $\modG$ has some countable base $\cU$.
	Recall that by assumption, $T$ is discrete in $\modG$.
	Therefore, we find for each~$t \in T$ a basic neighbourhood~$U_t \in \cU$ so that $U_t$ and $T \setminus \{t\}$ are disjoint;
	in particular, $U_t \neq U_{t'}$ for distinct $t, t' \in T$.
	Then $T$ is countable because $\cU \supseteq\left\{U_t \colon t \in T \right\}$ is countable.
	
	Since $T \cap \Omega(G)$ is countable, we may fix an enumeration $(\omega_n \colon n < \kappa)$ of $T \cap \Omega(G)$ where $\kappa:=\left| T \cap \Omega(G)\right|\leq \aleph_0$.
	Next, we recursively find for each $n < \kappa$ an $\omega_n$--$(T \setminus \{\omega_n\})$ cut $F_n$ of~$G$ such that the component $C_n$ of $G-F_n$ in which $\omega_n$ lives is disjoint from the component $C_m$ of $G-F_m$ in which $\omega_m$ lives for all $m < \kappa$ other than~$n$ (see also \cref{fig:MainProof} for a visualisation of the whole proof).
	
	Given any~$n<\kappa$, assume that we have already found suitable finite cuts $F_i$ for all $i < n$.
	Let $G_n$ be the graph obtained from~$G$ by contracting each component~$C_i$ for $i<n$ to a single vertex~$v_i$, keeping all the parallel edges that may arise.
	Since all cuts $F_i$ for $i<n$ are finite, the contraction minor $G_n$ is again locally finite. 
	Let
	\[T_n := \left(T \setminus \left\{\omega_i \colon i < n \right\}\right) \cup \left\{v_i \colon i < n \right\} \subseteq \hat{V}(G_n).\]
	Note that $\modGn = \modG\; / \left\{\overline{C_i} \colon i < n \right\}$, since all the $F_i$ are finite.
	
	Since $T$ is discrete in $\modG$ and the components $C_i$ for $i<n$ are disjoint, $T_n$ is discrete in $\modGn$.
	Therefore, the end $\omega_n$ is not contained in~$\overline{T_n \setminus \{\omega_n\}}$ where the closure is taken in $\modGn$.
	Thus, there is a smallest $\omega_n$--$(T_n\setminus \{\omega_n\})$ cut $F^*_n$ in $G_n$.
	This finite cut $F^*_n$ in $G_n$ defines a finite cut $F_n$ of $G$ of the same size.
	
	Finally, we observe that the component $C_n \subseteq G_n$ of $G_n-F^*_n$ in which $\omega_n$ lives does not contain any $v_i$ with $i < n$, since $F^*_n$ separates $\omega_n$ and $T_n \setminus\{\omega_n\} \supseteq \left\{v_i \colon i < n \right\}$.
	Thus, $C_n$ is also a component of $G- F_n$ and disjoint from each previous $C_i$.
	Altogether, this shows that the cut $F_n$ is as desired.
	
	Next, we simultaneously contract each component $C_n$ for $n < \kappa$ to a single vertex~$v_n$, again keeping all the parallel edges that may arise, and obtain a contraction minor $G_\kappa$ of $G$.
	As before, this contraction minor $G_\kappa$ is locally finite since all the cuts $F_n$ are finite. 
	Let
	\[T_\kappa := \left( T \setminus \Omega(G) \right) \cup \left\{v_n \colon n < \kappa \right\},\]
	and note that $T_\kappa \subseteq V(G_\kappa)$.
	Moreover, $T_\kappa$ is countable because it has the same size as $T$. 

    We show that there is no vertex set~$X\subset V(G)$ whose closure in~$\modG$ is disjoint from~$T$ and for which~$d_G(X)$ is an odd natural number:
    Let~$X$ be a set of vertices of~$G$ such that its closure in~$\modG$ is disjoint from~$T$ and~$d_G(X)$ is finite.
    It follows from the Jumping Arc Lemma~\cite{diestelbook}*{Proposition~8.6.3~(i)} that~$T$ lives in~$V(G)\setminus X$.
    Hence, $d_G(X) = |E_G(X,V(G)\setminus X)|$ is even by assumption.
 
	This yields that there is no vertex set $X\subset V(G_\kappa)$ whose closure in~$\vert G_\kappa\vert$ is disjoint from~$T_\kappa$ and for which $d_{G_\kappa}(X)$ is an odd natural number.
	Since $T_\kappa$ contains no ends, it further follows that there is no vertex set $X\subset V(G_\kappa)\setminus T_\kappa$ for which $d_{G_\kappa}(X)$ is an odd natural number.
	Therefore, we can apply \cref{t:LChInfinite} in~$G_\kappa$ to~$T_\kappa$ to obtain a set $\cP$ of edge-disjoint $T_\kappa$-paths in~$G_\kappa$ with the following property:
	For every vertex $t' \in T_\kappa$, there is a cut $F'_{t'}$ of~$G$ on the set of $t'$--$(T_\kappa\setminus\{t'\})$ paths in~$\cP$.
	Note that all cuts $F'_{t'}$ are finite, because $G_\kappa$ is locally finite and the paths in~$\cP$ are edge-disjoint.
	It remains to translate the set~$\cP$ of edge-disjoint $T_\kappa$-paths in~$G_\kappa$ into the desired set $\cA$ of pairwise edge-disjoint graphic $T$-arcs in~$\modG$.
	
	For every $n < \kappa$, the finite $\omega_n$--$(T \setminus \{\omega_n\})$ cut $F_n$ has smallest size.
	Let $S_n$ be the set of those endvertices of $F_n$ in $G - C_n$.
	By definition, $F_n$ is a minimal $S_n$--$\omega_n$ cut in $G$.
	So we can apply \cref{lemma:MengerVerticesEnd} to $S_n$ and $\omega_n$ to find a set $\cR_n$ of $\vert F_n\vert$ many edge-disjoint $S_n$--$\omega_n$ rays.
	Note that the rays in $\cR_n$ use only vertices of $V(C_n) \cup S_n$ since $F_n$ is a cut.
	In particular, there is for each edge in $F_n$ precisely one ray in $\cR_n$ which starts in this edge.
	
	Now every $T_\kappa$-path $P\in\cP$ in~$G_\kappa$ uniquely defines a graphic $T$-arc in~$\modG$ as follows: 
	If the first or last edge $e$ of $P$ runs between a contraction vertex $v_n$ and another vertex~$u$, then we replace it with the unique $u$--$\omega_n$ ray in $\cR_n$ which traverses the inner points of~$e$ and add the end~$\omega_n$ to~it.
	Here we allow one exception:
	if $P$ consists of just one edge $e$ between two vertices $v_n$ and $v_m$ in~$T_\kappa$, then we replace $e$ with the double ray which arises from the two rays in~$\cR_n$ and $\cR_m$ that start in~$e$, and add the ends $\omega_n$ and $\omega_m$ to it.
	In either case, let us write $A(P)$ for the graphic arc defined by~$P$ in this way.
	We claim that $\cA:=\{A(P)\mid P\in\cP\}$ is the desired set of graphic $T$-arcs.
	
	The arcs in~$\cA$ are edge-disjoint because the paths in~$\cP$ are edge-disjoint, the components $C_n$ are disjoint, and the rays in each set $\cR_n$ are edge-disjoint.
	It remains to show that for each $t \in T$, the number of arcs in~$\cA$ that link $t$ to $T\setminus\{t\}$ is equal to~$\lambda_G(t,T\setminus\{t\})$.
	Given any $t\in T$, let $t'\in T_\kappa$ be equal to $t$ if $t$ is a vertex, and let $t':=v_n$ if $t$ is an end~$\omega_n$.
	The finite $t'$--$(T_\kappa\setminus\{t'\})$ cut $F'_{t'}$ of~$G_\kappa$ witnesses that the set of $t'$--$(T_\kappa\setminus\{t'\})$ paths in~$\cP$ form a maximal sized edge-disjoint $t'$--$(T_\kappa\setminus\{t'\})$ path system in~$G_\kappa$.
	Since the $t$--$(T\setminus\{t\})$ arc system defined by~$\cA$ arises from this path system by replacing each path $P$ with the arc $A(P)$, the number of arcs in this system is equal to the size of the finite cut $F'_{t'}$.
	As the $t'$--$(T_\kappa\setminus\{t'\})$ cut $F'_{t'}$ of~$G_\kappa$ induces a $t$--$(T\setminus\{t\})$ cut of~$G$ of the same size, we have $\lambda_G(t,T\setminus\{t\})\le\vert F'_{t'}\vert$, which completes the proof.
\end{proof}

\begin{appendix}

\section{A degree condition for vertices and ends}\label{app:degrees}

We recall that an end $\omega$ of a locally finite graph~$G$ is \emph{even}, or \emph{has even degree}, if there is a finite vertex set $S\subset V(G)$ such that for every finite set $S'\supset S$ of vertices, the maximum number of edge-disjoint $S'$--$\omega$ rays is even;
otherwise, $\omega$ is \emph{odd} or \emph{has odd degree}.
We refer to \cite{bruhn2007end}*{§3} for a discussion of the parity of ends.
With the parity of ends at hand, it seems natural to extend the notion of \emph{inner-Eulerian} from finite to locally-finite infinite graphs and their ends, as follows.
A~locally finite graph~$G$ is \emph{inner-Eulerian} for a set $T\subset\hat{V}(G)$ if every vertex and every end in $\hat{V}(G)\setminus T$ have even degree in~$G$.

This raises the following three questions:
\begin{enumerate}[label=(\arabic*)]
    \item\label{Q1} Is it necessary to consider ends in this generalisation of `inner-Eulerian'?
    \item\label{Q2} Why does the premise of \cref{t:MainResult} not use this notion of `inner-Eulerian'?
    \item\label{Q3} How does the the notion of `inner-Eulerian' compare to the premise of \cref{t:MainResult}?
\end{enumerate}
The first question is answered by the following example.

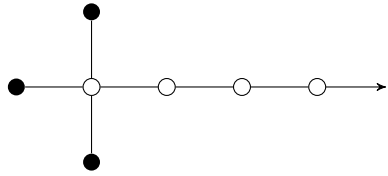
\begin{figure}[ht]
	\centering
	\begin{tikzpicture}[>=stealth',every node/.style={scale=0.62}]
		\node[circle, fill=black] (v1) at (0,0) {};
		\node[circle, draw] (v2) at (0,-1) {};
		\node[circle, fill=black] (v3) at (-1,-1) {};
		\node[circle, fill=black] (v4) at (0,-2) {};
		\node[circle, draw] (v5) at (1,-1) {};
		\node[circle, draw] (v6) at (2,-1) {};
		\node[circle, draw] (v7) at (3,-1) {};
		\node (end) at (4,-1) {};

		\draw[->] (v7) edge (end);
		
		\draw  (v1) edge (v2);
		\draw  (v2) edge (v3);
		\draw  (v2) edge (v4);
		\draw  (v2) edge (v5);
		\draw  (v5) edge (v6);
		\draw  (v6) edge (v7);
	\end{tikzpicture}
	\caption{$T$ consists of the black vertices} \label{fig:FailureTwo}
\end{figure}

\begin{example}
    \cref{t:MainResult} becomes wrong if we replace its premise on the vertex sets~$X$ with an `inner-Eulerian' condition that ignores end degrees.
    Indeed, let us consider the tree~$G$ in \cref{fig:FailureTwo}, and let $T$ consist of its leaves.
    All the non-leaves of~$G$ have even degree, but the unique end of~$G$ has degree one.
    Any set $\cA$ of pairwise edge-disjoint $T$-arcs has size at most one, so $\cA$ contains no $t$--$(T\setminus\{t\})$ arc for at least one $t\in T$, even though $\lambda_G(t,T\setminus\{t\})=1$.
\end{example}

We answer the second question twofold:
On the one hand, with the premise of \cref{t:MainResult} it is clear that \cref{t:MainResult} implies \cref{t:LChInfinite} for locally-finite graphs.
On the other hand,
the premise of $G$ being inner-Eulerian for~$T\subseteq\hat{V}(G)$ is more specific than the premise on the finite cuts in \cref{t:MainResult}, see \cref{lem:NoOddDegreeNoOddCl} below.
The latter also answers the third question.

\begin{lemma}[Infinite Handshaking Lemma \cite{bruhn2007end}*{Proposition 15}] \label{prop:InfiniteHandshake}
    The number of odd vertices and ends in a locally finite graph is even or infinite. 
\end{lemma}

\begin{lemma} \label{lem:NoOddDegreeNoOddCl}
	If $G$ is a locally finite graph which is inner-Eulerian for $T\subset\hat{V}(G)$, then
        every finite cut of $G$ such that $T$ lives on one of its sides is even.
\end{lemma}

\begin{proof}
    Assume for a contradiction that there is an odd finite cut $E_G(A,B)$ such that $T$ lives in~$A$.
    Consider the graph $H$ that arises from $G$ by contracting $G-B$ to a single vertex $v$, keeping parallel edges.
	Then $v$ has odd degree in~$H$, but no other vertex or end of $H$ has odd degree, which contradicts \cref{prop:InfiniteHandshake}. 
\end{proof}

The converse of \cref{lem:NoOddDegreeNoOddCl} fails if $G$ is a ray and~$T=V(G)$.
Hence \cref{t:MainResult} is more general than its corollary below (which follows with \cref{lem:NoOddDegreeNoOddCl}):

\begin{corollary} \label{t:MainResultinDegrees}
    Let $G$ be any locally-finite graph, and let $T\subseteq \hat{V}(G)$ be discrete in $\modG$ such that~$G$ is inner-Eulerian for~$T$.
    Then $\modG$ contains a set $\cA$ of pairwise edge-disjoint graphic $T$-arcs such that for every $t \in T$, the number of $t$--$(T\setminus\{t\})$ arcs in~$\cA$ is equal to~$\lambda_G(t,T\setminus\{t\})$.\qed
\end{corollary}

\begin{ack}
We thank the two referees for helpful comments and efficient reviews.
We are particularly grateful for the suggestion to reformulate the premise of \cref{t:MainResult} with cuts.
\end{ack}

\end{appendix}

\bibliographystyle{plain}
\bibliography{BIB}

\begin{bibdiv}
\begin{biblist}

\bib{aharoni2009menger}{article}{
      author={Aharoni, Ron},
      author={Berger, Eli},
       title={Menger’s theorem for infinite graphs},
        date={2009},
     journal={Inventiones mathematicae},
      volume={176},
      number={1},
       pages={1\ndash 62},
  note={\href{https://doi.org/10.1007/s00493-007-2149-0}{https://doi.org/10.1007/s00222-008-0157-3}},
}

\bib{bruhn2005menger}{article}{
      author={Bruhn, Henning},
      author={Diestel, Reinhard},
      author={Stein, Maya},
       title={Menger's theorem for infinite graphs with ends},
        date={2005},
     journal={Journal of Graph Theory},
      volume={50},
      number={3},
       pages={199\ndash 211},
  note={\href{https://doi.org/10.1002/jgt.20108}{https://doi.org/10.1002/jgt.20108}},
}

\bib{bruhn2007end}{article}{
      author={Bruhn, Henning},
      author={Stein, Maya},
       title={On end degrees and infinite cycles in locally finite graphs},
        date={2007},
     journal={Combinatorica},
      volume={27},
      number={3},
       pages={269},
  note={\href{https://doi.org/10.1007/s00493-007-2149-0}{https://doi.org/10.1007/s00493-007-2149-0}},
}

\bib{cherkassky1977asol}{article}{
      author={Cherkassky, Boris~V.},
       title={A solution of a problem on multicommodity ﬂows in a network
  ({E}nglish title), {R}eshenie odnoi zadachi o mnogoproduktovykh potokakh v
  seti ({R}ussian title)},
        date={1977},
     journal={Ekonomika i Matematicheskie Metody},
      volume={13},
      number={1},
       pages={143\ndash 151},
        note={(in Russian)},
}

\bib{diestel2003countable}{article}{
      author={Diestel, Reinhard},
       title={The countable {E}rd{\H{o}}s--{M}enger conjecture with ends},
        date={2003},
     journal={Journal of Combinatorial Theory, Series B},
      volume={87},
      number={1},
       pages={145\ndash 161},
  note={\href{https://doi.org/10.1007/s00493-007-2149-0}{https://doi.org/10.1016/s0095-8956(02)00032-1}},
}

\bib{diestel2010topapr2}{article}{
      author={Diestel, Reinhard},
       title={Locally finite graphs with ends: {A} topological approach, {II.
  A}pplications},
        date={2010},
     journal={Discrete Mathematics},
      volume={310},
      number={20},
       pages={2750\ndash 2765},
  note={\href{https://doi.org/10.1016/j.disc.2010.05.027}{https://doi.org/10.1016/j.disc.2010.05.027}},
}

\bib{diestel2011topapr1}{article}{
      author={Diestel, Reinhard},
       title={Locally finite graphs with ends: {A} topological approach, {I.
  B}asic theory},
        date={2011},
     journal={Discrete Mathematics},
      volume={311},
      number={15},
       pages={1423\ndash 1447},
  note={\href{https://doi.org/10.1016/j.disc.2010.05.023}{https://doi.org/10.1016/j.disc.2010.05.023}},
}

\bib{diestelbook}{book}{
      author={Diestel, Reinhard},
       title={Graph theory. 5th. vol. 173. graduate texts in mathematics},
   publisher={Springer},
        date={2016},
}

\bib{DiestelPersComm}{misc}{
      author={Diestel, Reinhard},
       title={Personal communication},
        date={2020},
}

\bib{diestel2004infinitecyc1}{article}{
      author={Diestel, Reinhard},
      author={Kühn, Daniela},
       title={On infinite cycles {I}},
        date={2004},
     journal={Combinatorica},
      volume={24},
      number={1},
       pages={69\ndash 89},
  note={\href{https://doi.org/10.1007/s00493-004-0005-z}{https://doi.org/10.1007/s00493-004-0005-z}},
}

\bib{diestel2004infinitecyc2}{article}{
      author={Diestel, Reinhard},
      author={Kühn, Daniela},
       title={On infinite cycles {II}},
        date={2004},
     journal={Combinatorica},
      volume={24},
      number={1},
       pages={91\ndash 116},
  note={\href{https://doi.org/10.1007/s00493-004-0005-z}{https://doi.org/10.1007/s00493-004-0006-y}},
}

\bib{diestel2004topological}{article}{
      author={Diestel, Reinhard},
      author={Kühn, Daniela},
       title={Topological paths, cycles and spanning trees in infinite graphs},
        date={2004},
     journal={European Journal of Combinatorics},
      volume={25},
      number={6},
       pages={835\ndash 862},
  note={\href{https://doi.org/10.1016/j.ejc.2003.01.002}{https://doi.org/10.1016/j.ejc.2003.01.002}},
}

\bib{diestel2012topapr3}{article}{
      author={Diestel, Reinhard},
      author={Sprüssel, Philipp},
       title={Locally finite graphs with ends: {A} topological approach, {III.
  F}undamental group and homology},
        date={2012},
     journal={Discrete Mathematics},
      volume={312},
      number={1},
       pages={21\ndash 29},
  note={\href{https://doi.org/10.1016/j.disc.2011.02.007}{https://doi.org/10.1016/j.disc.2011.02.007}},
}

\bib{lovcherinf}{article}{
      author={Jo{\'o}, Attila},
       title={The lov{\'a}sz-cherkassky theorem in countable graphs},
        date={2023},
     journal={Journal of Combinatorial Theory, Series B},
      volume={159},
       pages={1\ndash 19},
  note={\href{https://doi.org/10.1016/j.jctb.2022.11.001}{https://doi.org/10.1016/j.jctb.2022.11.001}},
}

\bib{lovasz1976some}{article}{
      author={Lov{\'a}sz, L{\'a}szl{\'o}},
       title={On some connectivity properties of {E}ulerian graphs},
        date={1976},
     journal={Acta Mathematica Academiae Scientiarum Hungarica},
      volume={28},
      number={1--2},
       pages={129\ndash 138},
  note={\href{https://doi.org/10.1007/bf01902503}{https://doi.org/10.1007/bf01902503}},
}

\end{biblist}
\end{bibdiv}
\end{document}